\documentclass[12pt, twoside, eqno]{article}
\usepackage{latexsym}
\usepackage{amssymb}
\usepackage{amsfonts}
\begin{document}
\begin{center}
{\bf \large On intra-regular and some left regular\\
$\Gamma$-semigroups }\bigskip

{\bf Niovi Kehayopulu, Michael Tsingelis }\end{center}\bigskip

\smallskip

\noindent{\bf Abstract} {\small We characterize the intra-regular 
$\Gamma$-semigroups and the left regular $\Gamma$-semigroups $M$ in 
which $x\Gamma M\subseteq M\Gamma x$ for every $x\in M$ in terms of 
filters and we prove, among others, that every intra-regular 
$\Gamma$-semigroup is decomposable into simple components, and every 
$\Gamma$-semigroup $M$ for which $x\Gamma M\subseteq M\Gamma x$ is 
left regular, is decomposable into left simple components.\medskip

\noindent {\bf AMS 2010 Subject Classification:} 20M99 
(06F99)\medskip

\noindent{\bf Keywords:} $\Gamma$-semigroup; left (right)
congruence; semilattice congruence; left (right) ideal; left (right) 
simple; simple; intra-regular; left (right) regular}

\section{Introduction and prerequisites}A structure theorem 
concerning the intra-regular semigroups, another one concerning some 
left regular semigroups have been given in [3]. These are the two 
theorems in [3]: {\bf II.4.9 Theorem.} The following conditions on a 
semigroup $S$ are equivalent: (1) Every $\cal N$-class of $S$ is 
simple; (2) Every ideal of $S$ is completely semiprime; (3) For every 
$x\in S$, $x\in Sx^2S$; (4) For every $x\in S$, $N(x)=\{y\in S \mid 
x\in SyS\}$; (5) ${\cal N}={\cal I}$; (6) Every ideal of $S$ is a 
union of $\cal N$-classes. {\bf II.4.5 Theorem.} The following 
conditions on a semigroup $S$ are equivalent: (1) Every $\cal 
N$-class of $S$ is left simple; (2) Every left ideal of $S$ is 
completely semiprime and two-sided; (3) For every $x\in S$, $x\in 
Sx^2$ and $xS\subseteq Sx$; (4) For every $x\in S$, $N(x)=\{y\in S 
\mid x\in Sy\}$; (5) ${\cal N}={\cal L}$; (6) Every left ideal of $S$ 
is a union of $\cal N$-classes.

Note that we always use the term ``semiprime" instead of ``completely 
semiprime" given by Petrich in [3]. So the condition (2) in the two 
theorems above should be read as ``Every ideal (resp. left ideal) of 
$S$ is semiprime", meaning that if $A$ is an ideal (resp. left ideal) 
of $S$, then for every $x\in S$ such that $x^2\in A$, we have $x\in 
A$. In the present paper we generalize these results in case of 
$\Gamma$-semigroups.

Let $M$ be a $\Gamma$-semigroup. An equivalence relation $\sigma$ on 
$M$ is called {\it left} (resp. {\it right}) {\it congruence} (on 
$M$) if $(a,b)\in\sigma$ implies $(c\gamma a, c\gamma b)\in\sigma$ 
(resp. $(a\gamma c, b\gamma c)\in\sigma$) for every $c\in M$ and 
every $\gamma\in\Gamma$. A relation $\sigma$
which is both left and right congruence on $M$ is called a
{\it congruence} on $M$. A congruence $\sigma$ on $M$ is called {\it 
semilattice congruence} if $(a\gamma b,b\gamma a)\in\sigma$ and 
$(a\gamma a,a)\in\sigma$ for every $a,b\in M$ and every $\gamma\in 
\Gamma$. A nonempty subset $A$ of $M$ is called a {\it left} (resp. 
{\it right}) {\it ideal} of $M$ if $M\Gamma A\subseteq A$ (resp. 
$A\Gamma M\subseteq A$). A subset $A$ of $M$ which is both a left and 
right ideal of $M$ is called an {\it ideal} of $M$. For an element 
$a$ of $M$, we denote by $L(a)$, $R(a)$, $I(a)$ the left ideal, right 
ideal and the ideal of $M$, respectively, generated by $a$, and we 
have $L(a)=a\cup M\Gamma a$, $R(a)=a\cup a\Gamma M$,
$I(a)=a\cup M\Gamma a\cup a\Gamma M\cup M\Gamma a\Gamma M$.
We denote by $\cal L$ the equivalence relation on $M$ defined by 
${\cal L}:=\{(a,b) \mid L(a)=L(b)\}$, by $\cal R$ the equivalence 
relation on $M$ defined by ${\cal R}:=\{(a,b) \mid R(a)=R(b)\}$ and 
by $\cal I$ the equivalence relation on $M$ defined by ${\cal 
I}:=\{(a,b) \mid I(a)=I(b)\}$. A nonempty subset $A$ of $M$ is called 
a {\it subsemigroup} of $M$ if $a,b\in A$ and $\gamma\in\Gamma$ 
implies $a\gamma b\in A$, that is, $A\Gamma A\subseteq A$. A 
subsemigroup $F$ of $M$ is called a {\it filter} of $M$ if $a,b\in F$ 
and $\gamma\in\Gamma$ such that  $a\gamma b\in F$ implies $a\in F$ 
and $b\in F$. We denote by $\cal N$ the relation on $M$ defined by 
${\cal N}:=\{(a,b) \mid N(a)=N(b)\}$ where $N(x)$ is the filter of 
$M$ generated by $x$ $(x\in M)$. It is well known that the relation 
$\cal N$ is a semilattice congruence on $M$. So, if $z\in M$ and 
$\gamma\in\Gamma$, then we have $(z\gamma z,z)\in {\cal N}$, 
$(z\gamma z\gamma z,z\gamma z)\in {\cal N}$, $(z\gamma z\gamma 
z\gamma z,z\gamma z\gamma z)\in {\cal N}$ and so on. A subset $A$ of 
$M$ is called {\it semiprime} if $a\in M$ and $\gamma\in\Gamma$ such 
that $a\gamma a\in A$ implies $a\in A$. A $\Gamma$-semigroup 
$(M,\Gamma,.)$ is called {\it left simple} if for every left ideal 
$L$ of $M$, we have $L=M$, that is, $M$ is the only left ideal of 
$M$. A subsemigroup $T$ of $M$ is called {\it left simple} if the 
$\Gamma$-semigroup $(T,\Gamma,.)$ (that is, the set $T$ with the same 
$\Gamma$ and the multiplication ``." on $M$) is left simple. Which 
means that for every left ideal $A$ of $T$, we have $A=T$. A 
subsemigroup of $M$ which is both left simple and right simple is 
called {\it simple}. If $M$ is a $\Gamma$-semigroup and $\sigma$ a 
semilattice congruence on $M$, then the class $(a)_{\sigma}$ of $M$ 
containing $a$ is a subsemigroup of $M$ for every $a\in M$. Let now 
$M$ be a $\Gamma$-semigroup and $\sigma$ a congruence on $M$. For 
$a,b\in M$ and $\gamma\in\Gamma$, we define $(a)_{\sigma}\gamma 
(b)_{\sigma}:=(a\gamma b)_{\sigma}$. Then the set 
$M/{\sigma}:=\{(a)_{\sigma} \mid a\in M\}$ is a $\Gamma$-semigroup as 
well. A $\Gamma$-semigroup $M$ is said to be a semilattice of simple 
semigroups if there exists a semilattice congruence $\sigma$ on $M$ 
such that the class $(x)_{\sigma}$ is a simple subsemigroup of $M$ 
for every $x\in M$.\section{Intra-regular $\Gamma$-semigroups}We 
characterize here the intra-regular $\Gamma$-semigroups in terms of 
filtres and we prove that every intra-regular $\Gamma$-semigroup is 
decomposable into simple subsemigroups.\medskip

\noindent{\bf Definition 1.} (cf. [2]) A $\Gamma$-semigroup $M$ is 
called
{\it intra-regular} if $$x\in M\Gamma x\gamma x\Gamma M$$ for
every $x\in M$ and every $\gamma\in \Gamma$. \medskip

\noindent{\bf Lemma 2.} (cf. [1]) {\it If M is a
$\Gamma$-semigroup, then $\cal I\subseteq \cal N$}.\medskip

\noindent{\bf Theorem 3.} {\it Let M be a $\Gamma$-semigroup. The
following are equivalent:

$(1)$ M is intra-regular.

$(2)$ $N(x)=\{y\in M \mid x\in M\Gamma y\Gamma M\}$ for every
$x\in M$.

$(3)$ $\cal N=\cal I$.

$(4)$ For every ideal I of M, we have $I=
\bigcup\limits_{x \in I} {(x)_{\cal N}}$.

$(5)$ $(x)_{\cal N}$ is a simple subsemigroup of M for
every $x\in M$.

$(6)$ $M$ is a semilattice of simple semigroups.

$(7)$ Every ideal of M is semiprime.\medskip

\noindent Proof.} $(1)\Longrightarrow (2)$. Let $x\in M$ and
$T:=\{y\in M \mid x\in M\Gamma y\Gamma M\}$. $T$ is a filter of
$M$. In fact: Take an element $\gamma\in\Gamma$
$(\Gamma\not=\emptyset)$. Since $M$ is intra-regular, we have
$$x\in M\Gamma x\gamma x\Gamma M=(M\Gamma x)\gamma x\Gamma 
M\subseteq
(M\Gamma M)\gamma x\Gamma M\subseteq M\Gamma x\Gamma M,$$then $x\in
T$, and $T$ is a nonempty subset of $M$. Let $a,b\in T$ and
$\gamma\in\Gamma$. Then $a\gamma b\in T$. Indeed: Since $b\in T$,
we have $x\in M\Gamma b\Gamma M$. Since $a\in T$, $x\in M\Gamma
a\Gamma M$. Since $M$ is intra-regular, we have
\begin{eqnarray*}x\in M\Gamma x\gamma x\Gamma M&\subseteq&M\Gamma 
(M\Gamma b\Gamma
M)\gamma (M\Gamma a\Gamma M)\Gamma M\\&=&(M\Gamma M)\Gamma (b\Gamma 
M\gamma M\Gamma a)\Gamma (M\Gamma M)\\&\subseteq&M\Gamma (b\Gamma
M\gamma M\Gamma a)\Gamma M.\end{eqnarray*}We prove that $b\Gamma
M\gamma M\Gamma a\subseteq M\Gamma (a\gamma b)\Gamma M$. Then we
have $$x\in M\Gamma {\Big (}M\Gamma (a\gamma b)\Gamma M{\Big )}\Gamma 
M\subseteq
M\Gamma (a\gamma b)\Gamma M,$$ and $a\gamma b\in T$. For this
purpose, let $b\delta u\gamma v\rho a\in b\Gamma M\gamma M\Gamma
a$, where $u,v\in M$ and $\delta, \rho\in\Gamma$. Since $M$ is
intra-regular, $b\delta u\gamma v\rho a\in M$ and
$\gamma\in\Gamma$, we have\begin{eqnarray*}b\delta u\gamma v\rho
a&\in&M\Gamma (b\delta u\gamma v\rho a)\gamma (b\delta u\gamma v\rho
a)\Gamma M\\&=&(M\Gamma b\delta u\gamma v)\rho (a\gamma b)\delta
(u\gamma v\rho a\Gamma M)\\&\subseteq&M\Gamma (a\gamma b)\Gamma
M,\end{eqnarray*}so $b\delta u\gamma v\rho a\in M\Gamma (a\gamma
b)\Gamma M$. Let $a,b\in M$ and $\gamma\in\Gamma$ such that $a\gamma 
b\in T$. Then $a,b\in T$. Indeed: Since $a\gamma b\in T$, we have

$x\in M\Gamma (a\gamma b)\Gamma M=M\Gamma a\gamma (b\Gamma 
M)\subseteq M\Gamma a\Gamma M\;\mbox { and }$

$x\in (M\Gamma a)\gamma b\Gamma M\subseteq M\Gamma b\Gamma M,$\\so 
$a,b\in T$.
Let now $F$ be a filter of $M$ such that $x\in F$.
Then $T\subseteq F$. Indeed: Let $a\in T$. Then $x\in M\Gamma
a\Gamma M$, so $x=u\gamma a\rho v$ for some $u,v\in M$,
$\gamma,\rho\in\Gamma$. Since $u, a\rho v\in M$, $u\gamma (a\rho
v)\in F$ and $F$ is a filter of $M$, we have $u\in F$ and $a\rho
v\in F$. Since $a,v\in M$, $a\rho v\in F$ and $F$ is a filter, we
have $a\in F$ and $v\in F$, so $a\in F$.\smallskip

\noindent$(2)\Longrightarrow (3)$. Let $(a,b)\in\cal N$. Then
$a\in N(a)=N(b)$. Since $a\in N(b)$, by (2), we have $b\in M\Gamma
a\Gamma M\subseteq a\cup M\Gamma a\cup a\Gamma M\cup M\Gamma
a\Gamma M=I(a)$. Since $I(a)$ is an ideal of $M$ containing $b$,
we have $I(b)\subseteq I(a)$. Since $b\in N(a)$, by symmetry, we
get $I(a)\subseteq I(b)$. Then $I(a)=I(b)$, and
$(a,b)\in\cal I$. Thus we have $\cal N\subseteq \cal I$. On the
other hand, by Lemma 2, $\cal I\subseteq\cal N$. Thus $\cal N=\cal
I$.\smallskip

\noindent$(3)\Longrightarrow (4)$. Let $I$ be an ideal of $M$. If
$y\in I$, then $y\in (y)_{\cal N}\subseteq \bigcup\limits_{x \in
I} {(x)_{\cal N}}$. Let $y\in (x)_{\cal N}$ for some $x\in I$. Then, 
by (3),
$(y,x)\in\cal N=\cal I$, so $I(y)=I(x)$. Since $x\in I$ and $I(x)$
is the ideal of $M$ generated by $x$, we have $I(x)\subseteq I$.
Thus we have $y\in I(y)=I(x)\subseteq I$, and $y\in I$.\smallskip

\noindent$(4)\Longrightarrow (5)$. Let $x\in M$. Since $\cal N$ is a 
semilattice
congruence on $M$, $(x)_{\cal N}$ is a subsemigroup of $M$. Let
$I$ be an ideal of $(x)_{\cal N}$. Then $I=(x)_{\cal N}$. In fact:
Let $y\in (x)_{\cal N}$. Take an element $z\in I$ and an element
$\gamma\in\Gamma$ $(I, \Gamma\not=\emptyset)$. The set $M\Gamma
z\gamma z\gamma z\Gamma M$ is an ideal of $M$. Indeed, it is a
nonempty subset of $M$, and we have

$M\Gamma (M\Gamma z\gamma z\gamma z\Gamma M)=(M\Gamma M)\Gamma 
z\gamma z\gamma z\Gamma M\subseteq M\Gamma z\gamma z\gamma z\Gamma M$ 
and

$(M\Gamma z\gamma z\gamma z\Gamma M)\Gamma M=M\Gamma z\gamma z\gamma 
z\Gamma (M\Gamma M)\subseteq M\Gamma z\gamma z\gamma z\Gamma M$.\\
By hypothesis, we have $M\Gamma z\gamma z\gamma z\Gamma M=
\bigcup\limits_{t \in M\Gamma z\gamma z\gamma z\Gamma M}
{(t)_{\cal N}}.$\\Since $z\gamma z\gamma z\gamma z\gamma z\in
M\Gamma z\gamma z\gamma z\Gamma M$, we have $(z\gamma z\gamma
z\gamma z\gamma z)_{\cal N}\subseteq M\Gamma z\gamma z\gamma
z\Gamma M$. Since $(z\gamma z,z)\in\cal N$ and $z\in I\subseteq
(x)_{\cal N}$, we have $(z\gamma z\gamma z\gamma z\gamma z)_{\cal
N}=(z)_{\cal N}=(x)_{\cal N}$. Then $y\in (x)_{\cal N}\subseteq
M\Gamma z\gamma z\gamma z\Gamma M$ and

$y=a\delta z\gamma z\gamma
z\xi b=(a\delta z)\gamma z\gamma (z\xi b)\; \mbox { for some }
a,b\in M, \;\delta,\xi\in\Gamma.$\\We prove that $a\delta z,\;z\xi
b\in (x)_{\cal N}$. Then, since $I$ is an ideal of $(x)_{\cal N}$,
we have $(a\delta z)\gamma z\gamma (z\xi b)\in (x)_{\cal N}\Gamma
I\Gamma (x)_{\cal N}\subseteq I$, and $y\in I$. We
have\begin{eqnarray*}a\delta z\in (a\delta z)_{\cal
N}&:=&(a)_{\cal N}\delta (z)_{\cal N}=(a)_{\cal N}\delta (y)_{\cal
N} \mbox { (since } (z)_{\cal N}=(x)_{\cal N}=(y)_{\cal
N})\\&=&(a)_{\cal N}\delta(a\delta z\gamma z\gamma z\xi b)_{\cal
N}\\&=&(a)_{\cal N}\delta(a)_{\cal N}\delta (z\gamma z\gamma z\xi
b)_{\cal N}\\&=&(a)_{\cal N}\delta (z\gamma z\gamma z\xi b)_{\cal
N}\mbox { (since } (a\delta a,a)\in {\cal N})\\&=&{\Big(}a\delta 
(z\gamma z\gamma z\xi b){\Big)}_{\cal N}\\&=&(y)_{\cal
N}=(x)_{\cal N}\end{eqnarray*}and\begin{eqnarray*}z\xi b\in (z\xi
b)_{\cal N}&:=&(z)_{\cal N}\xi (b)_{\cal N}=(y)_{\cal N}\xi
(b)_{\cal N}=(a\delta z\gamma z\gamma z\xi b)_{\cal N}\xi
(b)_{\cal N}\\&=&(a\delta z\gamma z\gamma z)_{\cal N}\xi (b)_{\cal
N}\xi (b)_{\cal N}\\&=&(a\delta z\gamma z\gamma z)_{\cal N}\xi
(b\xi b)_{\cal N}\\&=&(a\delta z\gamma z\gamma z)_{\cal N}\xi
(b)_{\cal N}\\&=&(a\delta z\gamma z\gamma z\xi b)_{\cal
N}=(y)_{\cal N}=(x)_{\cal N}.
\end{eqnarray*}
\noindent$(5)\Longrightarrow (6)$. Since $\cal N$ is a semilattice 
congruence on $M$. \smallskip

\noindent$(6)\Longrightarrow (7)$. Suppose $\sigma$ be a semilattice 
congruence on $M$ such that $(x)_\sigma$ is a simple subsemigroup of 
$M$ for every $x\in M$. Let $I$ be an ideal of $M$,
$x\in M$ and $\gamma\in\Gamma$ such that $x\gamma x\in I$. The set 
$I\cap (x)_{\sigma}$ is an ideal of
$(x)_\sigma$. In fact: Since $x\gamma x\in I$ and $x\gamma x\in 
(x)_\sigma$, the set $I\cap (x)_\sigma$ is a nonempty subset of 
$(x)_{\sigma}$ and, since $(x)_\sigma$ is a subsemigroup of $M$, we 
have

$(x)_\sigma\Gamma (I\cap (x)_{\sigma})\subseteq (x)_\sigma\Gamma 
I\cap (x)_\sigma\Gamma (x)_\sigma\subseteq M\Gamma
I\cap (x)_\sigma\subseteq I\cap (x)_\sigma$ and

$(I\cap (x)_\sigma)\Gamma (x)_\sigma\subseteq I\Gamma (x)_\sigma\cap 
(x)_\sigma\Gamma (x)_\sigma\subseteq I\Gamma M\cap
(x)_\sigma\subseteq I\cap (x)_\sigma.$\\Since $(x)_\sigma$ is
a simple subsemigroup of $M$, we have $I\cap (x)_\sigma=(x)_\sigma$, 
and $x\in I$.\smallskip

\noindent$(7)\Longrightarrow (1)$. Let $a\in M$ and
$\gamma\in\Gamma$. Then $a\in M\Gamma a\gamma a\Gamma M$. Indeed:
The set $M\Gamma a\gamma a\Gamma M$ is an ideal of $M$. This is
because it is a nonempty subset of $M$ and

$M\Gamma(M\Gamma a\gamma a\Gamma M)=(M\Gamma M)\Gamma a\gamma a\Gamma 
M\subseteq M\Gamma a\gamma a\Gamma M$,

$(M\Gamma a\gamma a\Gamma M)\Gamma M=M\Gamma a\gamma a\Gamma (M\Gamma 
M)\subseteq M\Gamma a\gamma a\Gamma M.$\\By hypothesis, $M\Gamma 
a\gamma a\Gamma M$ is semiprime. Since $(a\gamma a)\gamma (a\gamma 
a)\in M\Gamma a\gamma
a\Gamma M$, we have $a\gamma a\in M\Gamma a\gamma a\Gamma M$, and
$a\in M\Gamma a\gamma a\Gamma M$. Thus $M$ is intra-regular.
$\hfill\Box$\section{On some left regular $\Gamma$-semigroups}Again 
using filters, we characterize here the left regular 
$\Gamma$-semigroups $M$ in which $x\Gamma M\subseteq M\Gamma x$ for 
every $x\in M$ and we prove that this type of $\Gamma$-semigroups are 
decomposable into left simple components.
If $x\Gamma M\subseteq M\Gamma x$ for every $x\in M$, then $A\Gamma 
M\subseteq M\Gamma A$ for every $A\subseteq M$. Indeed: If $a\in A$, 
$\gamma\in\Gamma$ and $b\in M$, then $a\gamma b\in a\Gamma M\subseteq 
M\Gamma a\subseteq M\Gamma A$. Thus if $A$ is a left ideal of $M$, 
then $A$ is a right ideal of $M$ as well. As a consequence, the left 
regular $\Gamma$-semigroups in which $x\Gamma M\subseteq M\Gamma x$ 
for every $x\in M$, are left regular and left duo. We also remark 
that the left regular $\Gamma$-semigroups are intra-regular. Indeed: 
Let $a\in M$. Since $M$ is left regular, we have $a\in M\Gamma 
a\gamma a\subseteq M\Gamma (M\Gamma a\gamma a)\gamma a\subseteq 
M\Gamma a\gamma a\Gamma M.$ The right regular $\Gamma$-semigroups are 
also intra-regular, and the right regular $\Gamma$-semigroups for 
which $M\Gamma x\subseteq x\Gamma M$ for every $x\in M$ are right 
regular and right duo, and decomposable into right simple 
subsemigroups.\medskip

\noindent{\bf Definition 4.} (cf. [2]) A $\Gamma$-semigroup $M$ is 
called {\it left} (resp. {\it right}) {\it regular} if $x\in M\Gamma 
x\gamma x$ (resp. $x\in x\gamma x\Gamma M)$ for every $x\in M$ and 
every $\gamma\in\Gamma$.\medskip

\noindent{\bf Lemma 5.} (cf. [1]) {\it If M is a
$\Gamma$-semigroup, then $\cal L\subseteq \cal N$ and $\cal 
R\subseteq \cal N$}.\medskip

\noindent{\bf Theorem 6.} {\it Let M be a $\Gamma$-semigroup. The
following are equivalent:

$(1)$ M is left regular and $x\Gamma M\subseteq M\Gamma x$ for every 
$x\in M$.

$(2)$ $N(x)=\{y\in M \mid x\in M\Gamma y\}$ for every
$x\in M$.

$(3)$ $\cal N=\cal L$.

$(4)$ For every left ideal L of M, we have $L=
\bigcup\limits_{x \in L} {(x)_{\cal N}}$.

$(5)$ $(x)_{\cal N}$ is a left simple subsemigroup of M for
every $x\in M$.

$(6)$ M is a semilattice of left simple semigroups.

$(7)$ Every left ideal of M is semiprime and two-sided.}\medskip

\noindent{\it Proof.} $(1)\Longrightarrow (2)$. Let $x\in M$ and 
$T:=\{y\in M \mid x\in M\Gamma y\}$. The set $T$ is a filter of $M$ 
containing $x$. In fact: Take an element $\gamma\in\Gamma$ 
$(\Gamma\not=\emptyset)$. Since $M$ is left regular, we have$$x\in 
M\Gamma x\gamma x=(M\Gamma x)\gamma x\subseteq (M\Gamma M)\gamma 
x\subseteq M\Gamma x,$$then $x\in T$, and $T$ is a nonempty subset of 
$M$. Let $a,b\in T$ and $\gamma\in\Gamma$. Then $a\gamma b\in T$. 
Indeed: Since $b, a\in T$, we have $x\in M\Gamma b$ and $x\in M\Gamma 
a$. Since $M$ is left regular, we have\begin{eqnarray*}x\in M\Gamma 
x\gamma x&\subseteq& M\Gamma(M\Gamma b)\gamma (M\Gamma a)=(M\Gamma 
M)\Gamma (b\gamma M\Gamma a)\\&\subseteq&M\Gamma (b\gamma M\Gamma 
a).\end{eqnarray*}We prove that $b\gamma M\Gamma a\subseteq M\Gamma 
a\gamma b$. Then we have$$x\in M\Gamma (M\Gamma a\gamma b)=(M\Gamma 
M)\Gamma (a\gamma b)\subseteq M\Gamma (a\gamma b),$$and $a\gamma b\in 
T$. Let now $b\gamma u\mu a\in b\gamma M\Gamma a$ for some $u\in M$, 
$\mu\in\Gamma$. Since $M$ is left regular, we 
have\begin{eqnarray*}b\gamma u\mu a&\in&M\Gamma (b\gamma u\mu 
a)\gamma (b\gamma u\mu a)=(M\Gamma b\gamma u)\mu (a\gamma b)\gamma 
(u\mu a)\\&\subseteq& M\Gamma {\Big(}(a\gamma b)\Gamma 
M{\Big)}\\&\subseteq& M\Gamma (M\Gamma a\gamma b)\mbox { (since 
}x\Gamma M\subseteq M\Gamma x \;\forall x\in M)\\&\subseteq&M\Gamma 
a\gamma b.\end{eqnarray*}Let $a,b\in M$ and $\gamma\in\Gamma$ such 
that $a\gamma b\in T$. Then $a,b\in T$. Indeed: Since $a\gamma b\in 
T$, we have $x\in M\Gamma a\gamma b\subseteq (M\Gamma M)\Gamma 
b\subseteq M\Gamma b$, so $b\in T$. By hypothesis, $a\gamma b\in 
a\Gamma M\subseteq M\Gamma a$. Then
$x\in M\Gamma a\gamma b\subseteq M\Gamma (M\Gamma a)\subseteq M\Gamma 
a$, so $a\in T$. Let now $F$ be a filter of $M$ such that $x\in F$. 
Then $T\subseteq F$. Indeed: Let $a\in T$. Then $x\in M\Gamma a$, 
that is $x=u\rho a$ for some $u\in M$, $\rho\in\Gamma$. Since $u\in 
M$, $\rho\in\Gamma$, $u\rho a\in F$ and $F$ is a filter of $M$, we 
have $u\in F$ and $a\in F$, then $a\in F$.\smallskip

\noindent$(2)\Longrightarrow (3)$. Let $(a,b)\in\cal N$. Then $a\in 
N(a)=N(b)$. Since $a\in N(b)$, by (2), we have $b\in M\Gamma 
a\subseteq a\cup M\Gamma a=L(a)$, so $L(b)\subseteq L(a)$. Since 
$b\in N(a)$, by symmetry, we get $L(a)\subseteq L(b)$. Then we have 
$L(a)=L(b)$, and $(a,b)\in\cal L$. By Lemma 5, $\cal L\subseteq \cal 
N$, so $\cal L=\cal N$.\smallskip

\noindent$(3)\Longrightarrow (4)$. Let $L$ be a left ideal of $M$. If 
$y\in L$, then $y\in (y)_{\cal N}\subseteq \bigcup\limits_{x \in L} 
{{{(x)}_{\cal N}}}$. Let $y\in (x)_{\cal N}$ for some $x\in L$. Then, 
by (3), $(y,x)\in {\cal N}={\cal L}$, so $L(y)=L(x)$. Since $x\in L$ 
and $L(x)$ is the left ideal of $M$ generated by $x$, we have 
$L(x)\subseteq L$. Then $y\in L(y)=L(x)\subseteq L$, so $y\in 
L$.\smallskip

\noindent$(4)\Longrightarrow (5)$. Let $L$ be a left ideal of 
$(x)_{\cal N}$. Then $L=(x)_{\cal N}$. In fact: Let $y\in (x)_{\cal 
N}$. Take an element $z\in L$ and an element $\gamma\in \Gamma$ 
$(L,\Gamma\not=\emptyset$). Since $M\Gamma z\gamma z$ is a left ideal 
of $M$, by hypothesis, we have $M\Gamma z\gamma z=\bigcup\limits_{t 
\in M\Gamma z\gamma z} {(t)_{\cal N}}$. Since $z\gamma z\gamma z\in 
M\Gamma z\gamma z$, we have $(z\gamma z\gamma z)_{\cal N}\subseteq 
M\Gamma z\gamma z$. Since $(z\gamma z,z)\in{\cal N}$ and $z\in 
L\subseteq (x)_{\cal N}$, we have $(z\gamma z\gamma z)_{\cal 
N}=(z)_{\cal N}=(x)_{\cal N}$. Then $y\in (x)_{\cal N}\subseteq 
M\Gamma z\gamma z$, thus $y=a\mu z\gamma z$ for some $a\in M$ and 
$\mu\in\Gamma$. We prove that $a\mu z\in (x)_{\cal N}$. Then, since 
$L$ is a left ideal of $(x)_{\cal N}$, we have $(a\mu z)\gamma z\in 
(x)_{\cal N}\Gamma L\subseteq L$, and $y\in L$. We 
have\begin{eqnarray*}a\mu z\in (a\mu z)_{\cal N}&=&(a)_{\cal N}\mu 
(z)_{\cal N}=(a)_{\cal N}\mu (y)_{\cal N}\mbox { (since } (z)_{\cal 
N}=(x)_{\cal N}=(y)_{\cal N})\\&=&(a)_{\cal N}\mu(a\mu z\gamma 
z)_{\cal N}=(a)_{\cal N}\mu(a)_{\cal N}\mu (z\gamma z)_{\cal 
N}\\&=&(a)_{\cal N}\mu(z\gamma z)_{\cal N}=(a\mu z\gamma z)_{\cal 
N}\\&=&(y)_{\cal N}=(x)_{\cal N}.\end{eqnarray*}
$(5)\Longrightarrow (6)$. Since $\cal N$ is a semilattice congruence 
on $M$.\smallskip

\noindent $(6)\Longrightarrow (7)$. Let $\sigma$ be a semilattice 
congruence on $M$ such that $(x)_{\sigma}$ is a left simple 
subsemigroup of $M$ for every $x\in M$. Let $L$ be a left ideal of 
$M$ and $x\in M$, $\gamma\in\Gamma$ such that $x\gamma x\in L$. The 
set $L\cap (x)_{\sigma}$ is a left ideal of $(x)_{\sigma}$. Indeed:
The set $L\cap (x)_{\sigma}$ is a nonempty subset of $(x)_{\sigma}$ 
(since $x\gamma x\in L$ and $x\gamma x\in (x)_{\sigma})$ and
$$(x)_{\sigma}\Gamma (L\cap (x)_{\sigma})\subseteq (x)_{\sigma}\Gamma 
L\cap (x)_{\sigma} \Gamma (x)_{\sigma}\subseteq M\Gamma L\cap 
(x)_{\sigma}\subseteq L\cap (x)_{\sigma}.$$
Since $(x)_{\sigma}$ is a left simple subsemigroup of $M$, we have 
$L\cap (x)_{\sigma}=(x)_{\sigma}$, then $x\in L$. Thus $L$ is 
semiprime. Let now $L$ be a left ideal of $M$. Then $L\Gamma 
M\subseteq L$. Indeed: Let $y\in L$, $\gamma\in\Gamma$ and $x\in M$. 
Since $L$ is a left ideal of $M$, we have $x\gamma y\in M\Gamma 
L\subseteq L$. The set $L\cap (x\gamma y)_{\sigma}$ is a left ideal 
of $(x\gamma y)_{\sigma}$. Indeed:

$\emptyset\not=L\cap (x\gamma y)_{\sigma}\subseteq (x\gamma 
y)_{\sigma}\mbox { (since } x\gamma y\in L \mbox { and } x\gamma y\in 
(x\gamma y)_{\sigma})$ and
$$(x\gamma y)_{\sigma}\Gamma (L\cap (x\gamma y)_{\sigma})\subseteq 
(x\gamma y)_{\sigma}\Gamma L\cap (x\gamma y)_{\sigma}\Gamma (x\gamma 
y)_{\sigma}\subseteq M\Gamma L\cap (x\gamma y)_{\sigma}.$$
Since $(x\gamma y)_{\sigma}$ is left simple, we have $L\cap (x\gamma 
y)_{\sigma}=(x\gamma y)_{\sigma}=(y\gamma x)_{\sigma}$, so $y\gamma 
x\in L$.\smallskip

\noindent$(7)\Longrightarrow (1)$. Let $x\in M$ and 
$\gamma\in\Gamma$. Since $M\Gamma x\gamma x$ is a left ideal of $M$, 
by hypothesis it is semiprime. Since $(x\gamma x)\gamma (x\gamma 
x)\in M\Gamma x\gamma x,$ we have $x\gamma x\in M\Gamma x\gamma x$, 
and $x\in M\Gamma x\gamma x$, thus $M$ is left regular. Let now $x\in 
M$. Then $x\Gamma M\subseteq M\Gamma x$. Indeed: Since $M$ is left 
regular, we have $x\in M\Gamma x\gamma x\subseteq (M\Gamma M)\Gamma 
x\subseteq M\Gamma x$, so $M\Gamma x$ is a nonempty subset of $M$. In 
addition, $M\Gamma (M\Gamma x)=(M\Gamma M)\Gamma x\subseteq M\Gamma 
x$, so $M\Gamma x$ is a left ideal of $M$. By hypothesis, $M\Gamma x$ 
is a right ideal of $M$ as well. Since $M\Gamma x$ is an ideal of $M$ 
containing $x$, we have $I(x)\subseteq M\Gamma x$. On the other hand, 
$x\Gamma M\subseteq x\cup M\Gamma x\cup x\Gamma M\cup M\Gamma x\Gamma 
M=I(x)$. Thus we obtain $x\Gamma M\subseteq M\Gamma x$. $\hfill\Box$

\noindent The right analogue of Theorem 6 also holds, and we have the 
following:\medskip

\noindent{\bf Theorem 7.} {\it Let M be a $\Gamma$-semigroup. The
following are equivalent:

$(1)$ M is right regular and $M\Gamma x\subseteq x\Gamma M$ for every 
$x\in M$.

$(2)$ $N(x)=\{y\in M \mid x\in y\Gamma M\}$ for every
$x\in M$.

$(3)$ $\cal N=\cal R$.

$(4)$ For every right ideal R of M, we have $R=
\bigcup\limits_{x \in R} {(x)_{\cal N}}$.

$(5)$ $(x)_{\cal N}$ is a right simple subsemigroup of M for
every $x\in M$.

$(6)$ M is a semilattice of right simple semigroups.

$(7)$ Every right ideal of M is semiprime and two-sided.}
{\small
\bigskip

\noindent Quasigroups and Related Systems 23, no. 2 (2015), to 
appear

\end{document}